\documentclass[12pt]{amsart} 
\usepackage{amssymb}

\makeatletter

\usepackage{euscript}

\let\oldmarginpar\marginpar
\long\def\marginpar#1{\oldmarginpar{\tiny\raggedright#1\par}}
\@mparswitchfalse 
\def\page#1{\ignorespaces}

\def\cases#1{\left\{\,\vcenter{\normalbaselines\m@th
    \ialign{$##\hfil$&\quad##\hfil\crcr#1\crcr}}\right.}
\def\matrix#1{\null\,\vcenter{\normalbaselines\m@th
    \ialign{\hfil$##$\hfil&&\quad\hfil$##$\hfil\crcr
      \mathstrut\crcr\noalign{\kern-\baselineskip}
      #1\crcr\mathstrut\crcr\noalign{\kern-\baselineskip}}}\,}

\newtheorem{theorem}{Theorem}
\newtheorem{lemma}{Lemma}
\newtheorem{corollary}{Corollary}
\newtheorem{example}{Example}

\newtheorem{uncorollary}{Corollary} 
\newtheorem{untheorem}{Theorem} 

\theoremstyle{definition}

\newtheorem{remark}{Remark} 
 
\newtheorem{remarks}{Remarks}

\def\eqlabel#1{\label{eq#1}}
\def\eqref#1{(\ref{eq#1})}

\let\<\langle
\let\>\rangle

\def\tanh{\operatorname{tanh}}

\let\over\@@over
\let\atop\@@atop
\let\above\@@above
\let\overwithdelims\@@overwithdelims
\let\atopwithdelims\@@atopwithdelims
\let\abovewithdelims\@@abovewithdelims

\makeatother

\begin{document}

\title{A new variant of The Schwarz--Pick--Ahlfors Lemma}
\author{Robert Osserman}
\address{}
\curraddr{MSRI\\
1000 Centennial Drive\\
Berkeley, CA  94720-5070\\
USA}

\begin{abstract}
We prove a ``general shrinking lemma'' that resembles the
Schwarz--Pick--Ahlfors Lemma and its many generalizations, but differs
in applying to maps of a finite disk into a disk, rather than
requiring the domain of the map to be complete. The conclusion is that
distances to the origin are all shrunk, and by a limiting procedure we
can recover the original Ahlfors Lemma, that {\em  all} distances are
shrunk. The method of proof is also different in that it relates the
shrinking of the Schwarz--Pick--Ahlfors-type lemmas to the comparison
theorems of Riemannian geometry.
\end{abstract}

\thanks{The methods and results of this paper derive from a paper of
Antonio Ros \cite{R}, and in particular, from Lemma 6 of that paper.
Research at MSRI is supported in part by NSF grant DMS-9701755.}

\maketitle

We start by reviewing the history of Schwarz-type lemmas, with remarks
about the effects---some beneficial and some not---of successive
generalizations. 

There are minor variations in the way the Schwarz lemma is usually
stated. Here is one of the standard formulations.

\begin{lemma}[The Schwarz Lemma]
\label{1.l}
Let $ f(z) $ be analytic on a disk $ |z| < R_1 $ and suppose that $
|f(z)| < R_2 $ and $ f(0) = 0 $. Then
\begin{equation}
\eqlabel{1.e}
|f(z)| \leq \frac{R_2}{R_1}\, |z| \text{ for } |z| < R_1.
\end{equation}
\end{lemma}

It is also generally noted that strict inequality holds for every $ z
\neq 0 $ unless $f$ is of the special form
\begin{equation}
\eqlabel{2}
f(z) = \frac{R_2}{R_1} \, e^{i \alpha}z, \text{ for some real }
\alpha.
\end{equation}

As immediate corollaries, one has:

\begin{corollary}[Liouville's Theorem]
\label{1.c}
A bounded analytic function in the entire plane is constant.
\end{corollary}
\begin{proof}
$ R_2 $ is fixed, and $ R_1 $ may be chosen arbitrarily large.
\end{proof}

\begin{corollary}
\label{2.c}
If $ R_1 = R_2 $, then
\begin{equation}
\eqlabel{3.e}
|f'(0)| \leq 1 .
\end{equation}
\end{corollary}

A slightly less obvious, but still elementary corollary is

\begin{corollary}
\label{3}
If $R_1 = R_2$ and if $f$ maps the boundary to the boundary,
then at any point $b$ with $|b|=R_1$ where $f'(b)$ exists, one has
\begin{equation}
\eqlabel{4}
|f'(b)| \geq 1.
\end{equation}
\end{corollary}

The proof follows immediately from the fact that distances to the
origin are shrunk under $f$, and therefore distances from the
boundary are stretched. More precisely, for $t$ real, $0<t<1$, we
have $|f(tb)| \leq t|b|$, so that
$$
|f(tb)-f(b)| \geq R_1-tR_1 = |tb-b|
$$
from which (4) follows.

We will return later to the possible significance of this elementary
observation. Although we will not make use of it here, we note that
a refinement of the above argument gives a stronger and sharp boundary
equality; namely, with $R_1 = R_2 = 1$, if a single boundary point $b$ 
maps to the boundary and if $f'(b)$ exists, then 
$$
|f'(b)| \geq 1+\frac{1-|f'(0)|}{1+|f'(0)|}.
$$
(See Osserman \cite{O}.)
 
In 1916, Pick \cite{P} gave a new slant to Schwarz' Lemma that was
to have an enormous impact on future developments:

\begin{lemma}[Schwarz-Pick Lemma]
\label{2}
Let $f(z)$ be a holomorphic map of the unit disk D into the unit
disk. Then
\begin{equation}
\eqlabel{5}
\hat{\rho}(f(z_1),f(z_2))\leq \hat{\rho}(z_1,z_2) \mbox{ for all } z_1,z_2 \in D,
\end{equation}
where $\hat{\rho}$ refers to distances measured in the hyperbolic 
metric in $D$.
\end{lemma}

What Pick observed was that we can compose $f$ with linear fractional
transformations that are isometries of the hyperbolic plane, taking
$z_1$ to 0 and $f(z_1)$ to 0. Then (5) reduces to
\begin{equation}
\eqlabel{6}
\hat{\rho}(0,f(z_2)) \leq \hat{\rho}(0,z_2).
\end{equation}
But hyperbolic distance to the origin is a monotonic function of euclidean
distance, so that (6) is equivalent to (1) (with $R_1 =R_2=1$).

\begin{uncorollary}
If $\parallel \ \parallel _H$ denotes norm in the hyperbolic metric, then
\begin{equation}
\eqlabel{7}
\parallel df_z \parallel _H \leq 1 \mbox{ for all } z \in D
\end{equation}
and if $\gamma$ is any curve in $D$, then the length of the image of $\gamma$
under $f$ is less than or equal to the length of $\gamma$, both measured
in the hyperbolic metric.
\end{uncorollary}

For future reference let us note the explicit form of these
quantities. The hyperbolic metric is given by
\begin{equation}
\eqlabel{8}
d \hat{s}^2 = \left( \frac{2}{1 - |z|^2}\right)^2 \, |d z|^2
\end{equation}
and its Gauss curvature $ \hat{K} $ satisfies
\begin{equation}
\eqlabel{9}
\hat{K} \equiv -1.
\end{equation}
\page{7}
Integrating \eqref{8} yields 
\begin{equation}
\eqlabel{10}
\hat{\rho} (0,z) = \log \frac{1 + |z|}{1 - |z|} = 2 \tanh^{-1} \, |z|.
\end{equation}

Ahlfors' great insight \cite{A} was that the same conclusions would
hold far more generally.

\begin{lemma}[Schwarz--Pick--Ahlfors Lemma]
\label{3.l}
Let $f$ be a holomorphic map of the unit disk $D$ into a Riemann
surface $S$ endowed with a Riemannian metric $ d s^2 $ with Gauss
curvature $ K \leq  -1 $. Then the hyperbolic length of any curve in $D$
is at least equal to the length of its image. Equivalently,
\begin{equation}
\eqlabel{11}
\rho(f(z_1), f(z_2)) \leq \hat{\rho}(z_1, z_2) \quad \text{ for all }
z_1, z_2 \text{ in } D
\end{equation}
or
\begin{equation}
\eqlabel{12}
\|d f_z\| \leq 1 \text{ everywhere},
\end{equation}
where the norm is taken with respect to the hyperbolic metric on $D$
and the given metric on the image.
\end{lemma}

With these results as background, let us give examples of the kind of
finite versions one can prove. We recall that a {\em  geodesic disk of
radius $R$} on a surface is the diffeomorphic image of a euclidean
disk of radius $R$ under the exponential map. Equivalently, one has
geodesic polar coordinates:
\begin{equation}
\eqlabel{13}
d s^2 = d \rho^2 + G(\rho, \theta)^2 d \theta^2,
\end{equation}
where $ \rho $ represents distance to the center of the disk, and
\begin{equation}
\eqlabel{14}
G(0, \theta) = 0, \quad \frac{\partial G}{\partial \rho}(0, \theta) =
1, \quad G(\rho, \theta) >  0, \quad 0 < \rho < R.
\end{equation}

We shall use the following notation throughout this paper: 

{\bf Notation:} Let $f$ map the disk $ |z| < R $ into a geodesic disk
centered at $ f(0) $ on a surface $S$ with metric $ ds^2 $. Then
\begin{equation}
\eqlabel{15}
\rho(p) = \text{ distance on } S \text{ from } f(0) \text{ to } p
\end{equation}
\begin{equation}
\eqlabel{16}
\hat{\rho}(z) = \text{ distance from } 0 \text{ to } z \text{ with
respect to a metric } d \hat s^2 \text{ on } |z| < R.
\end{equation}

\begin{example} 
\label{1.ex}
Let $f$ be a holomorphic map of $ |z| < R_1 $ into a geodesic disk of
radius $ R_2 $ centered at $ f(0) $ on a surface $S$ with Gauss
curvature $ K \leq 0 $. Then
\begin{equation}
\eqlabel{17}
\rho (f(z)) \leq \frac{R_2}{R_1} \, |z|, \quad |z| < R_1.
\end{equation}
\end{example}

Note that this is a direct extension of the original Schwarz Lemma,
and it has exactly the same consequences:
\setcounter{corollary}{0}
\begin{corollary}
\label{1.corol}
Any holomorphic map of the entire plane into a geodesic disk on a
surface with $ K \leq 0 $ must be constant.
\end{corollary}
\begin{corollary}
\label{2.corol}
If $ R_2 \leq R_1 $, then $ \|d f_0\| \leq 1 $.
\end{corollary}
\begin{corollary}
\label{3.c}
If $ R_2 = R_1 $ and if at some point $ z $ with $ |z| = R_1 $, $
\rho(f(z)) = R_1 $ and $ d f_z $ exists, then 
\begin{equation}
\eqlabel{18}
\|d f_z\| \geq 1.
\end{equation}
\end{corollary}

\begin{remarks}
1. This example is a slightly more general form of the first part of
Lemma 6 of Ros \cite{R}; his proof goes through without change.

2. Cor. 1 is false for $ K >  0 $; stereographic projection is
a non-constant conformal map of the entire plane onto a geodesic disk
consisting of the sphere minus a point.
\page{11a}
\end{remarks}

\begin{example}
\label{2.ex}
Let $f$ map $ |z| < r < 1 $ into a geodesic disk of radius $ \rho_2 $
centered at $ f(0) $ on a surface $S$ whose Gauss curvature satisfies
$ K \leq -1 $. Let $ \rho_1 $ be the hyperbolic radius of $ |z| = r $;
i.e.,
$$
\rho_1 = \log \frac{1 + r}{1 - r}
$$
by {\rm \eqref{10}}. Then if $ \rho_2 \leq \rho_1 $ and $ d \hat{s} $
is the hyperbolic metric on $ |z| < 1 $,
\begin{equation}
\label{19}
\rho(f(z)) \leq \hat{\rho}(z) \quad \text{ for } |z| < r.
\end{equation}
\end{example}

\setcounter{corollary}{0}
\begin{corollary}
\label{1.coro}
Under the same hypotheses,
\begin{equation}
\eqlabel{20}
\|d f_0\| \leq 1.
\end{equation}
\end{corollary}

\begin{corollary}
\label{2.coro}
If furthermore, $ \rho_2 = \rho_1 $ and $f$ maps the boundary into the
boundary, then at any point $z$ on $ |z| = r $ where $ d f_z $ exists,
\begin{equation}
\eqlabel{21}
\|df_z\| \geq 1.
\end{equation}
\end{corollary}

Not that in both these examples we can only assert distance 
\page{12a}
shrinking from the center, unlike Schwarz--Pick and its descendants.
In fact, as \eqref{18} and \eqref{21} indicate, the reverse is likely
to be true near the boundary. However, the original Ahlfors version of
Schwarz--Pick turns out to be a consequence. 
\begin{corollary}[Schwarz--Pick--Ahlfors]
\label{3.coro}
{\rm (Lemma 3 above.)}
\end{corollary}

\begin{proof}
We show that \eqref{12} holds if $f$ maps the full disk $ |z| < 1 $
into $S$ with $ K \leq -1 $. If $S$ is not simply-connected we may
lift the map $f$ to a map $ \tilde{f} $ into the universal covering
surface $ \tilde{S} $ of $S$, in which case \eqref{12} is equivalent
to $ \|d \tilde{f}_z\| \leq 1 $ everywhere. So we may as well assume
that $S$ is simply-connected. 

Let $ z_1, z_2 $ be any two points in the unit disk. By composing $f$
with an isometry of the hyperbolic plane taking $ 0 $ to $ z_1 $, we
can assume that $ z_1 = 0 $, and \eqref{11} becomes
\begin{equation}
\eqlabel{22}
\rho(f(0), f(z)) \leq \hat{\rho}(0, z_2).
\end{equation}
Choose any $ r_0 $ such that
$$
|z_2| < r_0 < 1.
$$
Let 
$$
\rho_0 = \max_{|z| \leq r_0} \rho(f(z)).
$$
Since $S$ is simply connected and $ K < 0 $, there exist global
geodesic coordinates on the disk
$$
D_{\rho_0}: \rho(p) < \rho_0.
$$

We let $ \tilde{f}(\zeta) = f ( r_0 \zeta) :\{ |\zeta| < 1 \}\rightarrow
D_{\rho_0} $. Let $ d \tilde{s}$ be the hyperbolic metric in $ |\zeta|
< 1 $ and choose $ r_1 $, $|z_2| < r_1 < r_0 $ so that \page{14a} 
$$
\rho_1 = \tilde{\rho} \left(\frac{r_1}{r_0}\right) \geq \rho_0.
$$
(This is always possible, since $ \tilde{\rho}
\left(r_1/r_0\right) \rightarrow \infty $ as $ r_1 \rightarrow
r_0 $.) We may now apply our lemma to
$$
\tilde{f}: \{|\zeta| < \frac{r_1}{r_0}\} \rightarrow D_{\rho_0}
$$
to conclude $ \rho (\tilde{f}(\zeta)) \leq \tilde{\rho}(\zeta) $ for $
|\zeta| < r_1/r_0 $, and in particular for $ \zeta_2 = z_2/r_0
$. Then
$$
\rho(f(z_2)) = \rho(\tilde{f}(\zeta_2)) \leq \tilde{\rho}(\zeta_2).
$$
But $ r_0 $ can be chosen arbitrarily close to $ 1 $, and
$$
r_0 \rightarrow 1 \Rightarrow \zeta_2 \rightarrow z_2, \quad d
\tilde{s} \rightarrow d \hat{s}, \text{ and }
\tilde{\rho}(\zeta_2) \rightarrow \hat{\rho} (z_2).
$$
This proves \eqref{22} and therefore \eqref{11}, from which \eqref{12}
follows.
\end{proof}

The proof of Example \ref{2.ex} is an obvious analog of the proof of
Lemma 2.1 of \cite{OR}. But Example 2 is also a special case of
Theorem 1 below:~the ``general finite shrinking lemma.''

Before stating the general shrinking lemma, let us note some of the
generalizations of the Ahlfors Lemma \page{15a} that were made
subsequently. 

\begin{untheorem}
{\rm Yau (\cite{Y}, 1973).}  Let $ \hat{S}$ be complete, with $ \hat{K} \geq
-1 $, and let $f$ be a holomorphic map of $ \hat{S} $ into $S$, with $
K \leq -1 $.  Then $ \|d f_p\| \leq 1 $ for all $ p $ in $ \hat{S} $;
i.e., the length of every curve in $ \hat{S} $ is greater than or
equal to the length of its image.
\end{untheorem}

\begin{untheorem}
{\rm Troyanov (\cite{T}, 1991), Ratto, Rigoli and V\'eron (\cite{RRV}, 1994).}
Let $ \hat{S} $ be complete and let $f$ map $ \hat{S} $
holomorphically into $S$. Suppose that
\begin{equation}
\eqlabel{23}
K(f(p)) \leq \hat{K}(p), 
\end{equation}
\begin{equation}
\eqlabel{24}
K(f(p)) \leq 0,
\end{equation}
and that certain further restrictions hold on $K$, $ \hat{K} $, weaker than in 
Yau's theorem. Then $ \|df_p \| \leq 1 $ for all $ p $ in $
\hat{S} $. 
\end{untheorem}

We refer to the original papers for the exact \page{16a} hypotheses
in each case. What is of interest here is condition \eqref{23} which
represents the natural culmination of the line of investigation
initiated by Ahlfors. The underlying philosophy is that the more
negative the curvature, the more a holomorphic map will shrink
distances and curve lengths. Note that we are really comparing two
metrics on the same domain:~the original metric $ d \hat{s}^2 $
and the pullback of the metric $ d s^2 $ under $f$. In fact all
of the Ahlfors-type lemmas may be stated as comparison theorems
between two conformally related metrics, and again, the philosophy is
that the more negative the curvature, the shorter the curve lengths
in the metric.

This type of result seems oddly reminiscent, but in apparent reverse,
of the standard comparison theorems from Riemannian geometry, which
say roughly that the more negative the curvature the more certain
curves are stretched. Specifically, one \page{17a} has:

\begin{lemma}[Comparison Lemma]
Let $ d s^2 $ and $ d \hat{s}^2 $ be metrics given in
geodesic polar coordinates by 
$$
d s^2 = d \rho^2 + G(\rho, \theta)^2 d \theta^2
$$
$$
d \hat{s}^2 = d \rho^2 + \hat{G}(\rho, \theta)^2 d \theta^2.
$$
If
\begin{equation}
\eqlabel{25}
K(\rho, \theta) \leq \hat{K} (\rho, \theta), \quad 0 < \rho < \rho_0,
\end{equation}
then
\begin{equation}
\eqlabel{26}
\frac{1}{G} \, \frac{\partial G}{\partial \rho} \geq \frac{1}{\hat{G}}
\, \frac{\partial \hat{G}}{\partial \rho}
\end{equation}
and
\begin{equation}
\eqlabel{27}
G(\rho, \theta) \geq \hat{G} (\rho, \theta), \quad 0 < \rho < \rho_0.
\end{equation}
\end{lemma}
\noindent
Note that
\begin{equation}
\eqlabel{28}
G(\rho_1, \theta) = \frac{d s}{d \theta} \text{ along the geodesic
circle } \rho = \rho_1,
\end{equation}
so that \eqref{27}, \eqref{28} imply that
\begin{equation}
\eqlabel{29}
L (\rho_1) \geq \hat{L}(\rho_1), \quad 0 < \rho_1 < \rho_0
\end{equation}
where $ L (\rho), \hat{L}(\rho) $ refer to the length in their
respective metrics of \page{18a} geodesic circles of radius $ \rho $.

An obvious question is what relation, if any, exists between the
Ahlfors-type lemmas and the Riemannian comparison theorem. The answer
is two-fold; first, there is a heuristic argument, based on \eqref{18}
and \eqref{21}, which provides a link between the two, and second, we
can use the Riemannian comparison lemma to prove a general finite
shrinking lemma which contains our Example \ref{2.ex} above as a special
case, and therefore provides a new route to proving the original
Ahlfors Lemma.

Let \page{1b} us start with a brief look at the heuristic argument
relating the two forms of comparison. We have a geodesic disk $ \hat{D}
$ of radius $ \rho_1 $ on a surface with Riemannian metric
$$
d {\hat{s}}^2 = d \hat{\rho}^2 + \hat{G} (\hat{\rho},\theta)^2 d
\theta
,
$$
where for any point $ P $ in $ \hat{D} $, $ \hat{\rho}(P) = \text{
distance } $ between $ P $ and the center $O$ of the disk. We map $
\hat{D} $ conformally by $f$ into a surface $S$ with metric $ d s^2
$, and assume that the image lies in a geodesic disk $D$ of the same
radius centered at the point $ f(0) $. Under suitable curvature
restrictions we wish to show that \page{2b}
\begin{equation}
\eqlabel{30}
\rho (f(P) \leq \hat{\rho} (P), \quad \text{ for all } P \text{ in } D
\end{equation}
where $ \rho(Q) = \text{ distance} $ on $S$ from $ f(0) $ to $Q$. We
introduce geodesic polar coordinates
$$
d s^2 = d \rho^2 + G (\rho, \theta)^2 d \theta^2 , \quad 0 \leq \rho
< \rho_1, \quad 0 \leq \theta < 2 \pi
$$
on the image, and the curvature relation we assume is that
\begin{equation}
\eqlabel{31}
K(\rho, \theta) \leq \hat{K}(\hat{\rho}, \theta) \quad \text{ when }
\rho = \hat{\rho};
\end{equation}
that is, for each fixed $ \theta $, the curvature of the image
geodesic disk is at most equal to the curvature of the original at the
same distance from the center. Then what we want to show, inequality
\eqref{30}, is that each geodesic disk $ \hat{\rho} < c $, for $ c <
\rho_1 $, maps into the geodesic disk $ \rho < c $ in the image.
Heuristically, the image disk is likely to be largest when $f$
\page{3b} maps $ \hat{D} $ {\em onto} the full disk $D$. So let us
assume that $f$ is such a map, and $f$ takes the boundary, $
\hat{\rho} = \rho_1 $, to the boundary, $ \rho = \rho_1 $. Let us
further assume that $f$ is defined and conformal in a slightly larger
disk $ \hat{\rho} < \rho_0 $. Then the Riemannian comparison lemma
applies, and we have inequality \eqref{29}, which tell us that {\em
globally}, the map $f$ takes the geodesic circle $ \hat{\rho} = \rho_1
$ of length $ \hat{L} (\rho_1) $ onto a geodesic circle of greater or
equal length $ L(\rho_1) $; {\em locally}, by virtue of \eqref{28},
the inequality \eqref{27} tells us that under the map of $ \hat{\rho}
= \rho_1 $ to $ \rho = \rho_1 $ which relates points with \page{4b}
the same angular coordinate $ \theta $, we have
\begin{equation}
\eqlabel{32}
\frac{d s}{d \hat{s}} \geq 1.
\end{equation}
However, $f$ will not in general preserve $ \theta $, so that
inequality \eqref{29} tells us only that \eqref{32} holds {\em  on average}
where $ s $ and $ \hat{s} $ represent arclength along $ \rho
= \rho_1 $ and $ \hat{\rho} = \rho_1 $ {\em  under the map} $f$. The
final heuristic assumption is that \eqref{32} holds along the whole
curve $ \rho = \rho_1 $, under the map $f$. Then conformality of $f$
implies that the same inequality also holds in the radial direction,
so that along each ``radius'': $ \theta = \theta_0 $ of $ \hat{D} $,
we have
\begin{equation}
\eqlabel{33}
\left. \frac{d \rho}{d \hat{\rho}}\,\, \right|_{\hat{\rho} = \rho_1} \geq 1,
\end{equation}
where $ \rho(\hat{\rho}) $ is the function whose value is $ \rho(f(P))$ 
at \page{5b} the point $P$ in $ \hat{D} $ with coordinates $
(\hat{\rho} $, $\theta_0) $. Finally, we may, by a standard type of
argument dating at least to Ahlfors' original paper, assume that
we have {\em  strict} inequality in \eqref{27}, and therefore in
\eqref{32} and \eqref{33}, and then get weak inequality by going to
the limit. Then what \eqref{33} tells us is that points in $ \hat{D} $
near the boundary $ \hat{\rho} = \rho_1 $ move {\em  further} from the
boundary $ \rho = \rho_1 $ of $D$, so that they move {\em  closer} to
the center of $D$; in other words, \eqref{30} holds, in fact with
strict inequality, for points $P$ in some annular region near the
boundary of $ \hat{D} $. We are then back to our original situation on
a disk of smaller radius in $ \hat{D} $, and \page{6b} we may expect
the same kind of contraction \eqref{30} to extend.

In brief, then, the heuristic connection is that an equality like
\eqref{31} on Gauss curvature implies an {\em  expansion} of the
boundary $ \hat{\rho} = \rho_1 $, to $ \rho = \rho_1 $, which by
conformality of $f$ implies an expansion in the radial direction from
the boundary, or a movement of points {\em  toward} the center, and
therefore a {\em  contraction} in the sense of \eqref{30}.

We have not been able to turn this heuristic argument into a complete
proof under the full generality of \eqref{31}, but we can do so for a
very broad class of metrics, including those of Examples \ref{1.ex}
and \ref{2.ex};
\page{7b}
namely the case when $ d \hat{s}^2 $ has circular symmetry.

\begin{theorem}[General Finite Shrinking Lemma]
\label{1.t}
Let $ \hat{D} $ be a geodesic disk of radius $ \rho_1 $ with respect
to a metric $ d \hat{s}^2 $. Assume that $ d \hat{s}^2 $ is
circularly symmetric, so that
\begin{equation}
\eqlabel{34}
d \hat{s}^2 = d \hat{\rho}^2 + \hat{G}(\hat{\rho})^2 d
\theta^2, \quad 0 \leq \hat{\rho} < \rho_1,
\end{equation}
where $ \hat{G} $ depends on $ \hat{\rho} $ only, and not $ \theta $.

Let $f$ be a holomorphic map of $ \hat{D} $ into a geodesic disk $D$
of radius $ \rho_2 $ on a surface $S$, with center at the image under
$f$ of the center of $ \hat{D}$. If $ \rho_2 \leq \rho_1 $, and if
\begin{equation}
\eqlabel{35}
K(\rho, \theta) \leq \hat{K} (\hat{\rho}) \quad \text{ for } \rho =
\hat{\rho},
\end{equation}
then
\begin{equation}
\eqlabel{36}
\rho(f(P)) \leq \hat{\rho}(P) \quad \text{ for all } P \text{ in }
\hat{D}.
\end{equation}
\end{theorem}
\begin{proof}
First, let us assume, as we may, that the metric \page{8b} \eqref{34}
is represented as a conformal metric
\begin{equation}
\eqlabel{37}
d \hat{s}^2 = \hat{\lambda}(r)^2 |d z|^2, \quad |z| < R \leq
\infty,
\end{equation}
where
\begin{equation}
\eqlabel{38}
|dz|^2 = dr^2 + r^2 d \theta^2, \quad 0 \leq r < R
\end{equation}
is the euclidean metric on the disk. Comparing \eqref{34} with
\eqref{37}, \eqref{38}, we find
\begin{equation}
\eqlabel{39}
d \hat{\rho} = \hat{\lambda} (r) d r,
\end{equation}
and
\begin{equation}
\eqlabel{40}
\hat{G} (\hat{\rho}) = r \hat{\lambda} (r).
\end{equation}
Thus
\page{9b}
\begin{equation}
\eqlabel{41}
\hat{\rho} = h(r) =: \int^r_0 \hat{\lambda}(t)\, dt, \quad 0 \leq r < R,
\end{equation}
where $ h(r) $ is a monotone strictly increasing function, with
\begin{equation}
\eqlabel{42}
h(R) = \rho_1 .
\end{equation}
Hence $h$ has inverse
\begin{equation}
\eqlabel{43}
r = H(\hat{\rho}), \quad 0 \leq \hat{\rho} < \rho_1,
\end{equation}
also monotone strictly increasing, with
$$
H(\rho_1) = R.
$$
We next recall that in geodesic polar coordinates:
$$
d s^2 = d \rho^2 + G (\rho, \theta)^2 d \theta^2,
$$
the Gauss curvature $K$ is given by the formula
\begin{equation}
\eqlabel{44}
K(\rho, \theta) = - \frac{1}{G(\rho, \theta)} \,
\frac{\partial^2 G}{\partial \rho^2} (\rho, \theta),
\end{equation}
while the Laplacian $ \Delta $ with respect to the metric $ d s^2
$ of any function $ \varphi(\rho) $ is given by
\page{10b}
\begin{equation}
\eqlabel{45} 
\Delta \varphi = \frac{1}{G} \,
\left[\frac{\partial}{\partial \rho} (G \varphi ' (\rho))\right] =
\varphi '' (\rho) + \left(\frac{\partial}{\partial \rho} \log G\right)
\varphi '(\rho),
\end{equation}
and in particular, when $ \varphi(\rho) = \rho $,
\begin{equation}
\eqlabel{46}
\Delta \rho = \frac{\partial}{\partial \rho} \log G =
\frac{1}{G}\, \frac{\partial G}{\partial \rho};
\end{equation}
therefore, \eqref{45} can be written
\begin{equation}
\eqlabel{47}
\Delta \varphi = \varphi '' (\rho) + (\Delta \rho) \varphi ' (\rho).
\end{equation}
Comparing \eqref{25} and \eqref{26} with \eqref{35} and \eqref{46}, we
see that our hypotheses imply
\begin{equation}
\eqlabel{48}
\Delta \rho \mid_{\rho = c} \geq \hat{\Delta} \hat{\rho}
\mid_{\hat{\rho} = c} \quad \text{ for } 0 < \rho < \rho_2 .
\end{equation}
Since the function $ H $ in  \eqref{43} satisfies
$$
H'(\hat{\rho}) = \frac{dr}{d \hat{\rho}} = 1 / \frac{d \hat{\rho}}{d
r} >  0,
$$
it follows from \eqref{47} and \eqref{48} and the definition of $H$
that
\page{11b}
$$
\Delta \log H(\rho) \mid_{\rho = c} \, \geq \hat{\Delta} \log
H(\hat{\rho}) \mid_{\hat{\rho} = c} \, = \hat{\Delta} \log |z| = 0,
$$
since $ \log |z| $ is harmonic in the euclidean metric and therefore
in any conformal metric on $D$. Since $f$ is holomorphic, it follows
that the pullback of $ \log H (\rho) $ to $D$ also satisfies
$$
\Delta_z \log H(\rho(f(z))) \geq 0
$$
at all points where $ \rho(f(z)) \neq 0 $, while the function 
$$ 
\log H(\rho(f(z))) \rightarrow - \infty 
$$ 
as $ \rho(f(z)) \rightarrow 0 $. Let
$$
u(z) = \log \frac{H(\rho(f(z)))}{|z|} \quad \text{ for } 0 < |z| < R.
$$
Then $ \Delta_z u \geq 0 $ at all points where $ 0 < |z| < R $, $
\rho(f(z)) \neq 0 $. Also, $ u \rightarrow - \infty $ where $ \rho
(f(z)) \rightarrow 0 $. Hence $u$ is subharmonic for $ 0 < |z| < R $.
Furthermore near $ z = 0 $ if we represent
\page{12b}
the map $f$ by $ w = F(z) $ in terms of a local isothermal
parameter $w$ near $ f(0) $, with $ w = 0 $ at $ f(0) $, then
$$
H(\rho(f(z))) \sim \frac{\lambda(0)}{\tilde{\lambda}(0)} |F' (0) | |z|
$$
where $ d s^2 = \lambda^2(w)\, |d w|^2 $. Thus $ u(z) $ is bounded at
$ z = 0 $ if $ F'(0) \neq 0 $, and $ u(z) \rightarrow - \infty $ as $ z
\rightarrow 0 $ if $ F' (0) = 0 $. In either case, $ u(z) $ is
subharmonic in the full disk $ |z| < R $, and by the maximum principle,
$$
\log \frac{H(\rho(f(z)))}{|z|} = u(z) \leq \lim_{|z| \rightarrow R}
u(z) \leq \log \frac{H(\rho_2)}{R},
$$
or
\begin{equation}
\eqlabel{49}
H(\rho(f(z))) \leq \frac{H(\rho_2)}{H(\rho_1)}\, |z| \leq |z|
\end{equation}
since $ H(\rho_1) = R $, and $ \rho_2 \leq \rho_1 \Rightarrow H(\rho_2)
\leq H(\rho_1) $.

Thus
$$
\rho(f(z)) \leq h(|z|) = \hat{\rho} (z)
$$
which is \eqref{36}, and proves the Theorem.
\end{proof}

\begin{remark}
We have said that this theorem includes Examples \ref{1.ex} and
\ref{2.ex} as special cases, but although that is true of Example 2,
the theorem initially implies only the case of Example 1 where $ R_2
\leq R_1 $. However, it is easy to see, but somewhat awkward to state,
what happens when $ \rho_2 > \rho_1 $. We have
\end{remark}
\begin{theorem}[More General Finite Shrinking Lemma]
\label{2.t}
Under the same hypotheses and with the same notation as in the $ GFSL
$, except that $ \rho_2 \neq \rho_1 $, we have the inequality 
\begin{equation}
\eqlabel{50}
\rho(f(z)) \leq h \left( \frac{H(\rho_2)}{H(\rho_1)} \, |z| \right) \quad \text{
for } 0 \leq |z| < R,
\end{equation}
provided, in the case $ \rho_2 >  \rho_1 $, that we make the
additional \page{14b} assumption that the metric $ d \hat{s}^2 $
extends as a circularly symmetric metric to a larger disk $ |z| < R_2
$, with $ \hat{\rho}(R_2) = \rho_2 $, and that the inequality
{\rm \eqref{35}} holds whenever $ \rho = \hat{\rho} < \rho_2 $. 
\end{theorem}

The proof given for Theorem \ref{1.t} goes through unchanged till the
first inequality in {\rm \eqref{49}}, which is equivalent to {\rm
\eqref{50}}.

In the case of Example 1, where $ d
\hat{s}^2 $ is the euclidean metric, $ h(r) = r $, $ H(r) = r $,
and \eqref{50} reduces to \eqref{17}.


\begin{thebibliography}{}

\bibitem [A]{A}
Ahlfors, L.V.``An extension of Schwarz's Lemma.'' {\em Trans.
Amer. Math. Soc.}{\bf 43} (1938) 359-364.

\bibitem [O]{O} Osserman, R. ``A sharp Schwarz inequality on
the boundary.'' (MSRI preprint No. 1997-110)

\bibitem [OR]{OR} Osserman, R. and Ru, M.  ``An estimate for the
Gauss curvature of minimal surfaces in $\bf{R}^m$ whose Gauss map
omits a set of hyperplanes.'' {\em J. Differential Geometry}
{\bf 46} (1997) 578-593.

\bibitem [P]{P} Pick, G. ``\"Uber eine Eigenschaft der konformen
Abbildung kreisf\"ormiger Bereiche." {\em Math. Annalen} {\bf 77}
(1916).

\bibitem [R]{R} Ros, A. ``The Gauss map of minimal surfaces.''
(preprint)

\bibitem [RRV]{RRV} Ratto, R., Rigoli, M. and V\'eron, L. ``Conformal
immersions of complete Riemannian manifolds and extensions of the
Schwarz lemma.'' {\em Duke Math. J.} {\bf 74} (1994) 223-236.

\bibitem [T]{T} Troyanov, M. ``The Schwarz Lemma for nonpositively
curved Riemannian surfaces.'' {\em Manuscripta Math.} {\bf 72} (1991)

\bibitem [Y]{Y} Yau, S.-T. ``Remarks on conformal transformations''
{\em J. Diff. Geometry} {\bf 8} (1973) 369-381.

\end{thebibliography}
\end{document}